# Image compression predicated on recurrent iterated function systems

Chol-Hui Yun[*], Metzler W.[a] and Barski M. [a]

[*] *Faculty of Mathematics, **Kim Il Sung** University, Pyongyang, D. P. R. Korea*
[a] *Faculty of Mathematics, University of Kassel, Kassel, F. R. Germany*

May, 2008

**Abstract**. Recurrent iterated function systems (RIFSs) are improvements of iterated function systems (IFSs) using elements of the theory of Marcovian stochastic processes which can produce more natural looking images. We construct new RIFSs consisting substantially of a vertical contraction factor function and nonlinear transformations. These RIFSs are applied to image compression.

## 1. Introduction

Fractal image compression (FIC) was introduced by Barnsley and Sloan (1988, [2]) and after that has been widely studied by many scientists. FIC is based on the idea that any image (for example that of human face) contains selfsimilarities, that is, it consists of small parts similar to itself or to some big part in it. So in FIC iterated function systems are used for modeling. Jacquin (1992, [12]) presented a more flexible method of FIC than Barnsley's, which is based on recurrent iterated function systems (RIFSs) introduced first by him. By his method, the original image is partitioned into small non-overlapping blocks called regions (or ranges) and big overlapping blocks called domains. For every region, a similar domain and some transformation are searched. Fisher (1994,[10]) improved the partition of Jacquin. Bouboulis etc.(2006, [7]) introduced an image compression scheme using fractal interpolation surfaces which are attractors of some RIFSs. RIFSs which have been used in image compression schemes consist of transformations which have a constant vertical contraction factor. This factor has a great influence on the shape of attractor. It is obvious that each point of a region has a different contraction ratio. So, replacing the

---

vertical contraction constant by a contraction function can give a more flexible construction so that the number of domains and regions needed for image coding can be reduced and the quality of the decoded image can be better, that is, the FIC will be improved.

We present a new construction of RIFSs which have vertical contraction factor functions and their application to image compression.

## 2. Recurrent bivariate IFSs on rectangular grids

The idea of RIFSs was introduced by Jacquin ([12]). A recurrent iterated function system (RIFS) is defined as a pair of a collection $w_1, \cdots, w_N$ of Lipschitz mappings in a complete metric space (i.e. an IFS) and a row irreducible stochastic matrix $P = (p_{ij})_{N \times N}$ satisfying $\sum_{j=1}^{N} p_{ij} = 1$ for all $i \in \{1, \cdots, N\}$ and additionally for any $i, j \in \{1, \cdots, N\}$, there exist $i_1, \cdots, i_k$ such that $p_{i i_1} \cdot p_{i_1 i_2} \cdot \cdots \cdot p_{i_k j} > 0$ (Barnsley, Elton and Hardin [5]). The attractor A of a RIFS is computed as follows: By a stochastic matrix $P$, for an initial $k_0 \in \{1, \cdots, N\}$ a sequence $\{k_i\}_0^\infty$ is given, which obeys $p_{ab} = \Pr\{k_{i+1} = b | k_i = a\}$. With this, we get a sequence of transformations $\{w_{k_i}\}_0^\infty$ that generates an orbit $\{q_i\}_0^\infty$ for an initial $q_0 \in \mathbf{R}^3$, that is, $k_{i+1}$ is chosen with the probability that $k_{i+1} = j$ equals $p_{k_i j}$ and we have $q_{i+1} = w_{k_{i+1}}(q_i)$. In the matrix $(p_{kl})$, $p_{kl}$ gives the possibility of applying the transformation $w_l$ to the point in state $k$, so that the system transits to state $l$.

The attractor A is defined as the limit set of these orbits, which consists of the points whose every neighbourhood contains infinitely many $q_i$ for almost all orbits. By Barnsley, Elton and Hardin ([5]) the existence, uniqueness and characterization of this limit set A was proved. We present a new construction of a RIFS with a given data set on a grid.

We construct a local IFS $\{\mathbf{R}^3; w_{ij} = (L_{ij}, F_{ij}), i = 1, \cdots, m, j = 1, \cdots, n\}$ over the grid, whose general definition is as follows ([3]).

**Definition 1** Let $(\mathbf{X}, d)$ be a compact metric space. Let $R$ be a nonempty subset of $\mathbf{X}$. Let $w: R \to \mathbf{X}$ and let s be a real number with $0 \leq s < 1$. If

$$d(w(x), w(y)) \leq s \cdot d(x, y) \quad \text{for all } x, y \text{ in } R,$$

then $w$ is called *a local contraction mapping* on $(\mathbf{X}, d)$. The number $s$ is called a contractivity factor for $w$.

**Definition 2** Let $(\mathbf{X}, d)$ be a compact metric space, and let $w_i: R_i \to \mathbf{X}$ be a local contraction mapping on $(\mathbf{X}, d)$, with contractivity factor $s_i$ for $i = 1, 2, \cdots, N$, where $N$ is a finite positive integer. Then

$$\{w_i: R_i \to \mathbf{X} : i = 1, 2, \cdots, N\}$$

is called *a local iterated function system (local IFS)* (or *partitioned IFS* ([10])). The number $s = \max\{s_i: i = 1, 2, \cdots, N\}$ is called the contractivity factor of the local IFS.

Let the data set on the rectangular grid be

$$\{(x_i, y_j, z_{ij}) \in \mathbf{R}^3 ; i = 0, 1, \cdots, m, j = 0, 1, \cdots, n\}$$

such that $x_0 < x_1 < \cdots < x_m$, $y_0 < y_1 < \cdots < y_n$. Let denote

$N_{mn} = \{1,\cdots,m\} \times \{1,\cdots,n\}$, $I_x = [x_0, x_m]$, $I_y = [y_0, y_n]$,

$I_{x_i} = [x_{i-1}, x_i]$, $I_{y_j} = [y_{j-1}, y_j]$, $E = I_x \times I_y$

$E_{ij} = I_{x_i} \times I_{y_j}$ (which we call the *region*), for $(i,j) \in N_{mn}$,

We choose big rectangulars (called *domains*) $\widetilde{E}_k = \widetilde{I}_{x_k} \times \widetilde{I}_{y_k}$, where $\widetilde{I}_{x_k} = [x'_{k_{k-1}}, x'_{k_k}]$ ($\widetilde{I}_{x_k} = \cup I_{x_\alpha}$), $\widetilde{I}_{y_k} = [y'_{k_{k-1}}, y'_{k_k}]$ ($\widetilde{I}_{y_k} = \cup I_{y_\beta}$) and $x_i - x_{i-1} < x'_{k_k} - x'_{k_{k-1}}$, $y_j - y_{j-1} < y'_{k_k} - y'_{k_{k-1}}$, for $k = 1, \cdots, l$. Then, there exist $x_{k_1}, x_{k_2} \in [0, m]$ ($y_{l_1}, y_{l_2} \in [0, n]$) such that $x'_{k_{k-1}} = x_{k_1}$ and $x'_{k_k} = x_{k_2}$ ($y'_{k_{k-1}} = y_{l_1}$ and $y'_{k_k} = y_{l_2}$), which we denote by $\gamma_x(x'_{k_{k-1}})$ and $\gamma_x(x'_{k_k})$ ($\gamma_y(y'_{k_{k-1}})$ and $\gamma_y(y'_{k_k})$) respectively.

We define the contraction transformations $L_{ij}: \widetilde{E}_k \to E_{ij}$ for $(i,j) \in N_{mn}$ by

$$L_{ij}(x, y) = \left(L_{x_i}(x), L_{y_j}(y)\right),$$

where $L_{x_i}: \widetilde{I}_{x_k} \to I_{x_i}$, $L_{y_j}: \widetilde{I}_{y_k} \to I_{y_j}$ are contraction transformations with contractivity factors $a_{x_i}$, $a_{y_j}$ obeying for any $(i,j) \in N_{mn}$

$$L_{x_i}: \{x'_{k_{k-1}}, x'_{k_k}\} \to \{x_{i-1}, x_i\},$$

$$L_{y_j}: \{y'_{k_{k-1}}, y'_{k_k}\} \to \{y_{j-1}, y_j\}.$$

The functions $F_{ij}: \widetilde{E}_k \times \mathbf{R} \to \mathbf{R}$, for $(i,j) \in N_{mn}$ are defined as follows:

$$F_{ij}(x, y, z) = d_{ij}\left(L_{ij}(x, y)\right)\left(z - g_k(x, y)\right) + h_{ij}\left(L_{ij}(x, y)\right),$$

where $|d_{ij}(x, y)| < 1$, $g_k$, $h_{ij}$ are continuous Lipschitz mappings on $\widetilde{E}_k$, $E_{ij}$ with the Lipschitz constants $L_{g_k}$, $L_{h_{ij}}$ respectively satisfying

$$g_k\left(x'_{k_\alpha}, y'_{k_\beta}\right) = z_{\sigma\left(x'_{k_\alpha}, y'_{k_\beta}\right)}, \quad \text{for} \quad \alpha, \beta \in \{k-1, k\},$$

$$h_{ij}(x_a, y_b) = z_{ab}, \quad \text{for} \quad (a, b) \in \{i-1, i\} \times \{j-1, j\},$$

where $\sigma\left(x'_{k_\alpha}, y'_{k_\beta}\right) = \sigma\left(\gamma_x(x'_{k_\alpha}), \gamma_y(y'_{k_\beta})\right) = \sigma(x_s, y_t) = (s, t)$, that is, $g_k$ goes through 4 endpoints of the domain $\widetilde{E}_k$, $h_{ij}$ goes through 4 endpoints of the region $E_{ij}$.

Then, the functions $F_{ij}$ satisfy 'join-up' conditions, for $\alpha, \beta \in \{k-1, k\}$,

$$F_{ij}\left(x'_{k_\alpha}, y'_{k_\beta}, z_{\sigma\left(x'_{k_\alpha}, y'_{k_\beta}\right)}\right) = z_{st},$$

where $L_{x_i}(x'_{k_\alpha}) = x_s$, $L_{y_j}(y'_{k_\beta}) = y_t$, $(s, t) \in \{i-1, i\} \times \{j-1, j\}$.

Thus, any $w_{ij}$ maps the endpoints of the $k$-th domain to the endpoints of the $(i,j)$-th region. These transformations $w_{ij}$ are contractions for all $(i,j) \in N_{mn}$

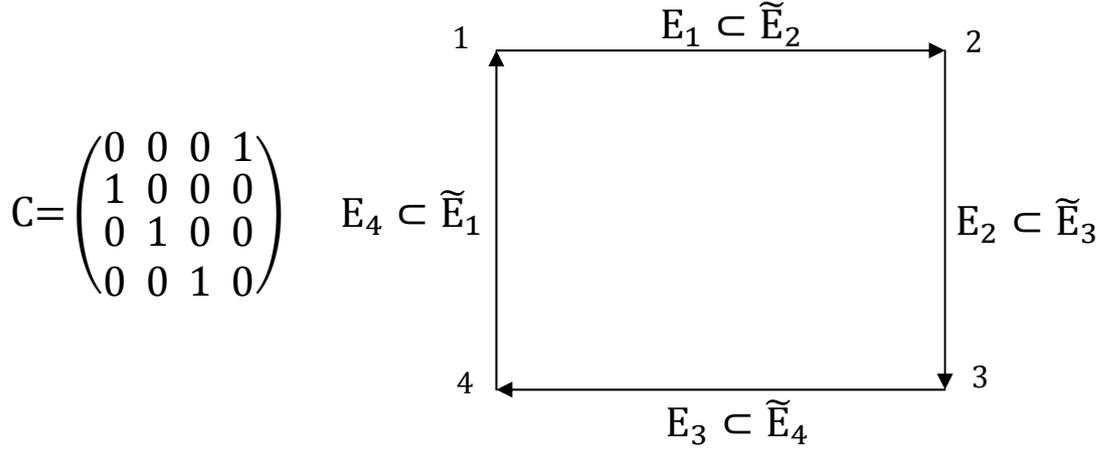
**Figure 1** *Connection matrix and it's directed graph.*

with respect to some metric $\rho$ which is equivalent to the Euclidean metric on $\mathbf{R}^3$. This metric $\rho$ is defined on $\mathbf{R}^3$ for $(x, y, z)$, $(x', y', z') \in \mathbf{R}^3$ by

$$\rho\big((x, y, z), (x', y', z')\big) = |x - x'| + |y - y'| + \theta |z - z'|,$$

where

$$\theta = \frac{1 - \text{Max}\{a_{x_i},\ a_{y_j};\ i = 1, \cdots, m,\ j = 1, \cdots, n\}}{2 L_F}$$

and $L_F = \text{Max}\{L_{F_{ij}};\ i = 1, \cdots, m,\ j = 1, \cdots, n\}$, where $L_{F_{ij}}$ are Lipschitz constants of the functions $F_{ij}$. The contractivity of the transformation $w_{ij}$ is $\text{Max}\{a, L_F\}$, where

$$a = \frac{1 + \text{Max}\{a_{x_i},\ a_{y_j};\ i = 1, \cdots, m,\ j = 1, \cdots, n\}}{2}$$

We enumerate the set $N_{mn} = \{1, \cdots, m\} \times \{1, \cdots, n\}$ by a injective mapping $\tau: \{1, \cdots, m\} \times \{1, \cdots, n\} \to \{1, \cdots, m \cdot n\}$ and denote $M = \tau(N_{mn})$. For simplicity, we denote $(i, j) = \tau^{-1}(k)$ by $k$ and $m \cdot n$ by $N$. The following is an example of a connection matrix $C = (c_{kl})_{N \times N}$ which can be defined as follows:

$$c_{kl} = \begin{cases} 1 & \text{if } E_l \subseteq \widetilde{E}_k \\ 0 & \text{otherwise} \end{cases}$$

The connection matrix $C$ shows that the transformation $w_k$ can follow the transformation $w_l$ iff $c_{kl}$ (See Fig. 1).

Then, we have a row irreducible stochastic matrix $P = (p_{kl})_{N \times N}$ given by:

$$c_{kl} = \begin{cases} 1 & \text{iff } p_{kl} > 0 \\ 0 & \text{iff } p_{kl} = 0 \end{cases},$$

Thus, the RIFS is given as a pair consisting of the above local IFS and the row irreducible stochastic matrix $P$. The existence, uniqueness and characteri-zation of the attractor A of this RIFS has been proved by Barnsley Elton and Hardin ([5]).

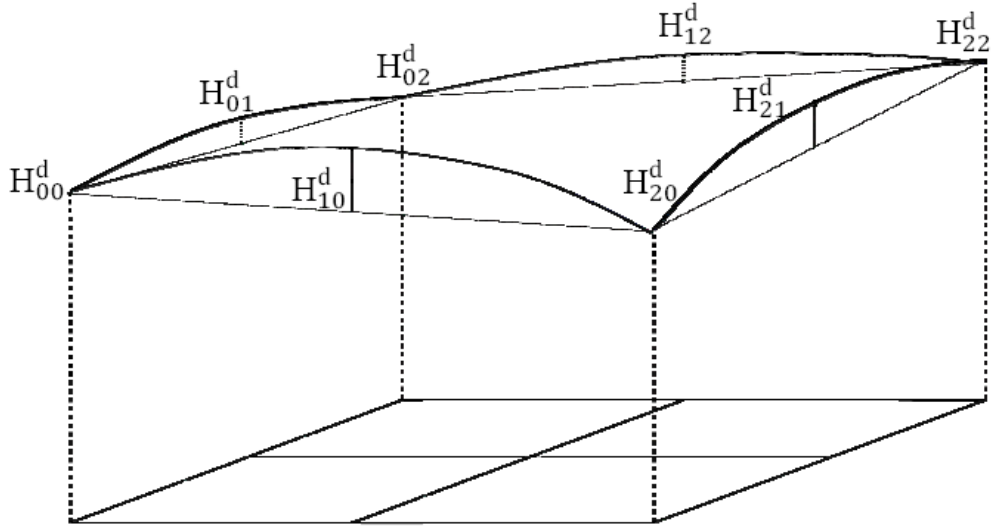

**Figure 2** $H_{ij}^d$ on the domain.

## 3. Application: An algorithm for image compression using RIFSs

Let us denote an image value at the position $(x,y) \in E$ by $I(x,y)$.

3.1 *Construction of the vertical contractivity factor functions*

In this section, we present a method of constructing the vertical contractivity factor function $d_{ij}(x,y)$ on the region $E_{ij}$ by the image values on the region $E_{ij}$ and on the corresponding domain $\widetilde{E}_k$ $((i,j) \in N_{mn}, k \in \{1, \cdots, l\})$.

The vertices $(x'_{k_{k-1}}, y'_{k_{k-1}})$, $(x'_{k_{k-1}}, y'_{k_k})$, $(x'_{k_k}, y'_{k_{k-1}})$, $(x'_{k_k}, y'_{k_k})$ of the domain $\widetilde{E}_k$ are mapped to any one of the region $E_{ij}$, $(x_{i-1}, y_{j-1})$, $(x_{i-1}, y_j)$, $(x_i, y_{j-1})$, $(x_i, y_j)$ by the contraction transformation $L_{ij} : \widetilde{E}_k \to E_{ij}$. Choose $y = y'_{k_{k-1}}$. For $x = x'_{k_{k-1}}, x'_{k_{k-1}} + \varepsilon_1^x, x'_{k_{k-1}} + \varepsilon_2^x, \cdots, x'_{k_{k-1}} + \varepsilon_{\delta_x}^x = x'_{k_k}$ ($\varepsilon_a^x \in \mathbf{R}^+$, $a = 1, \cdots, \delta_x$, $\delta_x \in \mathbf{N}$), we calculate the absolute vertical distance between the image value $I(x,y)$ and the straight line through the endpoints $(x'_{k_{k-1}}, y'_{k_{k-1}}, I(x'_{k_{k-1}}, y'_{k_{k-1}}))$ and $(x'_{k_k}, y'_{k_{k-1}}, I(x'_{k_k}, y'_{k_{k-1}}))$, and denote them by $H_{0,0}^d, H_{1,0}^d, \cdots, H_{\delta_x,0}^d$ respectively. In general, choose $y = y'_{k_{k-1}} + q$ for $q = 0, \varepsilon_1^y, \cdots, \varepsilon_{\delta_y}^y$, $y'_{k_{k-1}} + \varepsilon_{\delta_y}^y = y'_{k_k}$, $\varepsilon_b^y \in \mathbf{R}^+$, $b = 1, \cdots, \delta_y$, $\delta_y \in \mathbf{N}$, $H_{1,q}^d, \cdots, H_{\delta_x-1,q}^d$ are computed similarly.

Now, choose $x = x'_{k_{k-1}}$. For $y = y'_{k_{k-1}}, y'_{k_{k-1}} + \varepsilon_1^y, \cdots, y'_{k_{k-1}} + \varepsilon_{\delta_y}^y$, analogously calculate the absolute vertical distance between the image value

$I(x, y)$ and the straight line through the endpoints $\left(x'_{k_{k-1}}, y'_{k_{k-1}}, I(x'_{k_{k-1}}, y'_{k_{k-1}})\right)$ and $\left(x'_{k_{k-1}}, y'_{k_k}, I(x'_{k_{k-1}}, y'_{k_{k-1}})\right)$ and denote all these distances by $H^d_{0,0}, H^d_{0,1}, \cdots, H^d_{0,\delta_y}$ respectively. By the same method, $H^d_{\delta_x,0}, H^d_{\delta_x,1}, \cdots, H^d_{\delta_x,\delta_y}$ are computed. In the case where , $H^d_{k,l} = 0$ (excepting $H^d_{0,0}$, $H^d_{0,\delta_y}$, $H^d_{\delta_x,0}$, $H^d_{\delta_x,\delta_y}$), changing the value $\varepsilon^x_a$ or $\varepsilon^y_b$, we calculate that one again. Then, we get the set $H^d = \{H^d_{k,l}: k = 0, 1, \cdots, \delta_x,\ l = 0, 1, \cdots, \delta_y\}$ (See Fig. 2).

Similarly, in the region $E_{ij}$, for any point $(x, y)$, where $x = L_{x_i}(x')$, $y = L_{y_j}(y')$, $x' = x'_{k_{k-1}}$, $x'_{k_{k-1}} + \varepsilon^x_1, \cdots, x'_{k_{k-1}} + \varepsilon^x_{\delta_x} = x'_{k_k}$, $y' = y'_{k_{k-1}}$, $y'_{k_{k-1}} + \varepsilon^y_1$, $\cdots, y'_{k_{k-1}} + \varepsilon^y_{\delta_y} = y'_{k_k}$, the set $H^r = \{H^r_{i,j}: i = 0, 1, \cdots, \delta_x,\ j = 0, 1, \cdots, \delta_y\}$ is constructed.

With the sets $H^d$, $H^r$, we get a set
$$D_{ij} = \{d_{ij} = H^r_{i,j}/H^d_{i,j} : i = 0, 1, \cdots, \delta_x, j = 0, 1, \cdots, \delta_y\},$$
where $d_{00} = d_{0\delta_y} = d_{\delta_x 0} = d_{\delta_x \delta_y} = 0$.

We construct the contractivity factor function $d(x, y)$ by an interpolation function for the data set $D_{ij}$.

### 3.2 *Encoding algorithm*

Image encoding is a process to get reduced informations needed to preserve the data of an image. We introduce a way for getting the reduced information by a RIFS. We partition the original image into non-overlapping blocks on a rectangular grid and construct a grid of non-overlapping domains which are $a^2 (a \geq 2$ is integer) times as large as the regions . For each region, we look for the most similar domain to that as follows: By the above algorithm we calculate the contractivity factor function $d(x, y)$ on the region and the interpolation functions $g$ on the domain and $h$ on the region, respectively. Next we compute the distance between the image values on the region and the values mapped from the image values on the domain by a contraction transformation with regard to some metric. If the distance is less than a given tolerance, then we save the number of the domain and the data set for the contractivity factor function. If there is no suitable domain, then we divide equally the region into smaller rectangular (new regions), that is, repartition and repeat the above procedure for each smaller region. Finally, we get the reduced information, the endpoints of the regions, the number of the most similar domain for each region and the data set for the contractivity factor function.

### 3.3 *Decoding algorithm*

We calculate the interpolation functions $d_{ij}(x, y)$, $g_{ij}(x, y)$, $h_{ij}(x, y)$ from the information of the encoded image and compute the transformations $F_{ij}(x, y, z)$ out of them, that is, we determine the IFS.

The decoding algorithm is based on the Deterministic Iteration Algorithm (DIA) as in Barnsley [1] and it is similar to that in Bouboulis's [7], Section 3.3. The only difference compared with [7] is the following simplification: For any region, we find out if its center is already drawn. If yes, then we go to the next region. If no, then we search for an appropriate transformation of the IFS to draw this center and additionally the four midpoints of the edges of the region, and then we go to the next region.

**Table 1** *Some results obtained for Lena by the scheme in this paper.*

| Compression Time(CT)(second) | PSNR(dB) | Compression Ratio(CR) |
|---|---|---|
| 32 | 47.6 | 11.7 : 1 |
| 24 | 40.0 | 27.4 : 1 |
| 13 | 46.3 | 11.7 : 1 |
| 9 | 38.0 | 27.4 : 1 |
| 4 | 43.8 | 11.7 : 1 |
| 3 | 35.5 | 27.4 : 1 |

## 4. Experiments

We use a $513 \times 513 \times 8$ Lena image as original image in our experiment. According to the above algorithm, while this image is encoded and decoded, we measure the time needed to compress using a computer with a 1.86 GHz CPU clock and Windows XP. The experimental results show that our method can compress the image very fast, holding high quality of the decoded image (See Table 1). In Fig. 3-7, some images decoded by our scheme are shown.

## 5. Results & Discussion

In this paper we present a refined method of constructing an RIFS and its application to the fractal image compression (FIC). Its decoding scheme is simplified compared with earlier FIC schemes. Encoding is improved by a more complicated interpolation method based on the idea of using a vertical contractivity factor function instead of a constant contractivity factor.

The performance of FIC is significantly improved using our schemes. Now the compression ratio (CR) is not higher. To improve CR we are planning to add some kind of lossless compression and entropy coding compression.

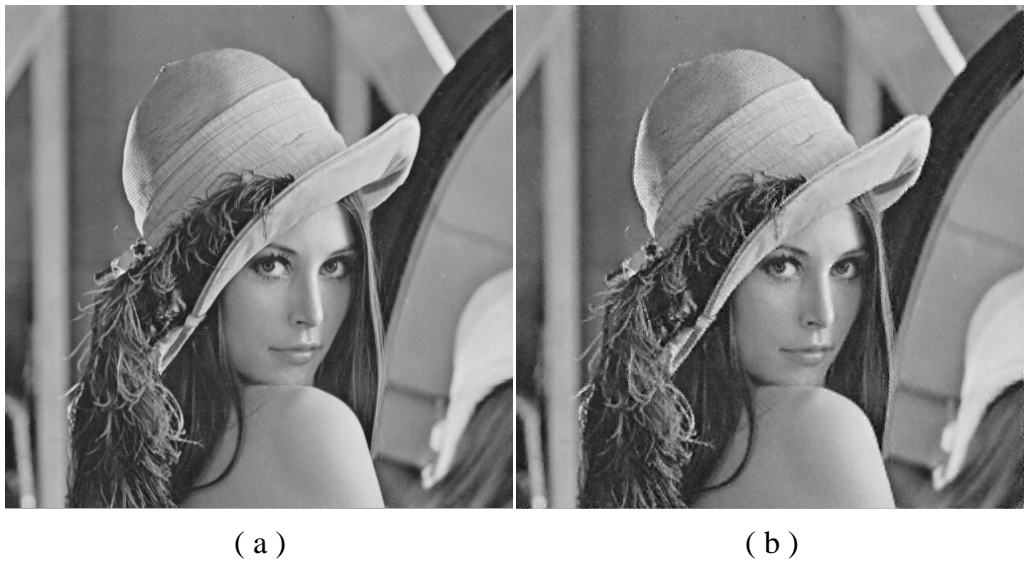

( a )          ( b )

**Figure 3** *( a ) The 513×513 original image of Lena. (b): The decoded image of Lena using our scheme at CT=32 seconds, PSNR=47.6dB, CR=11.7:1.*

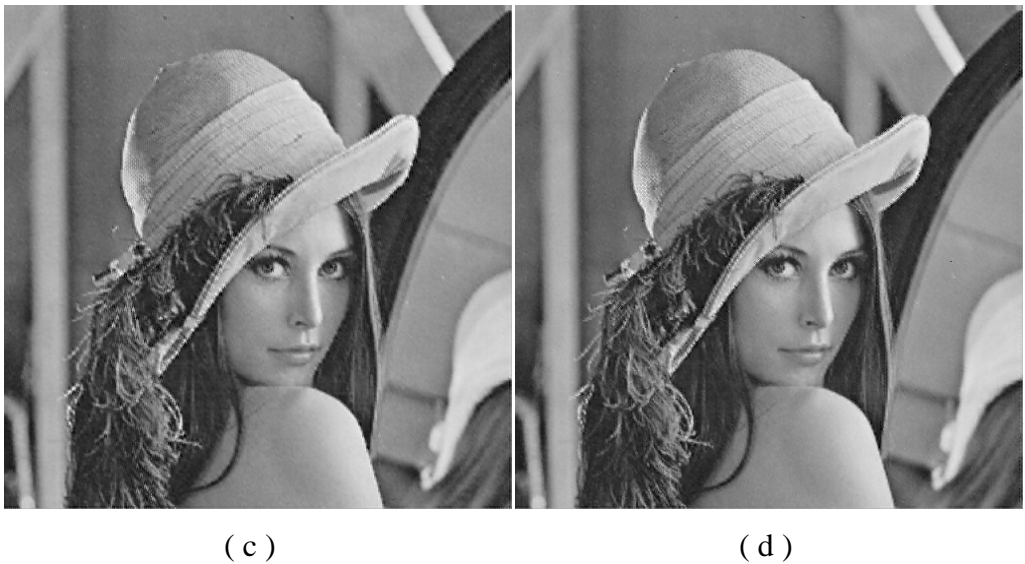

( c )          ( d )

**Figure 4** *(c): The decoded image of Lena at CT=4 seconds, PSNR=43.8, CR=11.7:1. (d): The decoded image of Lena at CT=13 seconds, PSNR=46.3 dB, CR=11.7:1.*

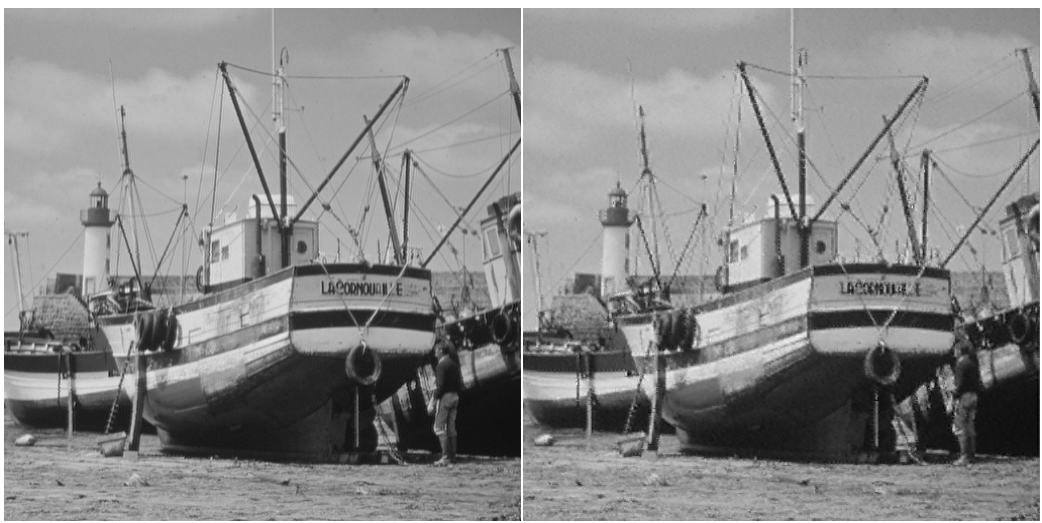

( a )     ( b )

**Figure 5** *( a ) : The 513×513 original image of the Boat. (b): The decoded image at CT=34 seconds, PSNR=43.6dB, CR=11.7:1*

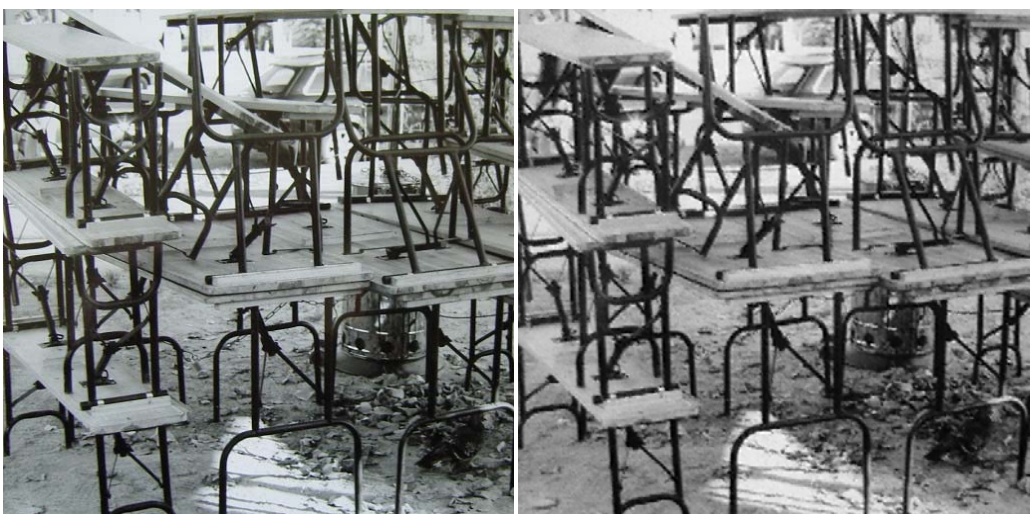

( a )     ( b )

**Figure 6** *( a ) : The 1025 ×1025 original image of the piled up benches. (b):The decoded image at CT=33 seconds, PSNR=36.5dB, CR=18.75:1.*

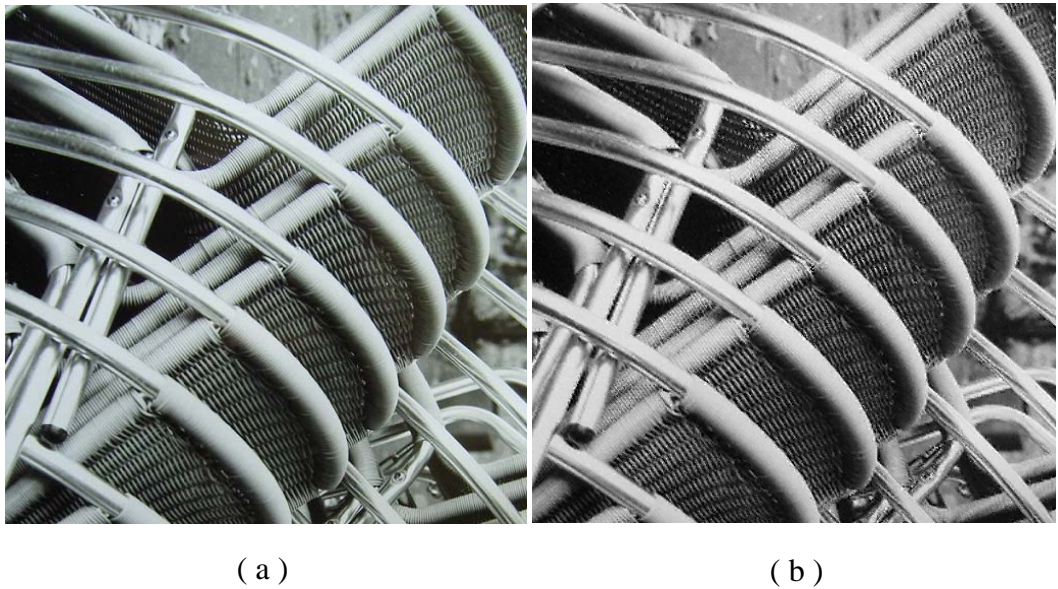

( a )                 ( b )

**Figure 7** *( a ): The 513×513 original image of the tower of chairs. (b): The decoded image at CT=11 seconds, PSNR=36.3dB, CR=11.7:1.*